\documentclass{article}
\usepackage[russian,english]{babel}
\usepackage{amssymb,amsmath,amsfonts}
\usepackage{graphics,dcolumn,bm,epic,eepic,float}
\usepackage{graphicx} 
\usepackage{epstopdf} 
\usepackage{cite}

\usepackage{bm} 
\usepackage{subfigure}
\usepackage{color} 

\usepackage{algorithm2e}
\usepackage[pdftex]{hyperref}

\usepackage{color, colortbl}
\definecolor{Yellow}{rgb}{1,1,0}
\definecolor{Grey}{rgb}{.87,.87,.87}
\definecolor{Purple}{rgb}{.8,.0,1.0}

\usepackage[english]{babel}
\begin{document}
\title{Painlev\'e Classification Of Polynomial 
Ordinary Differential Equations Of Arbitrary
Order And Second Degree}

\date{}

\author{Stanislav Sobolevsky
	\thanks{SENSEable City Laboratory, Massachusetts Institute of Technology, 77 Massachusetts Avenue 9-221, Cambridge, MA 02139, USA}
	\thanks{Institute of Lifelong Education of Belarusian State University, 15 Moskovskaya str, Minsk, Belarus}
	\thanks{E-mail: stanly@mit.edu}
}

\maketitle

\begin{abstract}
The problem of Painlev\'e classification of ordinary differential equations lasting since the end of XIX century saw significant advances for the limited equation order, however not that much for the equations of higher orders. In this work we propose the complete Painlev\'e classification for ordinary differential equations of the arbitrary order 
with right-hand side being a quadratic form on the dependent variable and all of its derivatives. The total of seven classes of
the equations with Painlev\'e property have been found. Five of them having the order up to four are already known. Sixth one of the other up to five also appears to be integrable in the known functions. While the only seventh class of the unrestricted order appears to be linearizable.
The classification employs a novel general necessary condition for the Painlev\'e property proven in the paper, potentially having a broader application for the Painlev\'e classification of other types of ordinary differential equations. 
\end{abstract}

\section*{Introduction}
The Painleve\'e classification is one of the long-lasting problems of analytic theory of differential equations rooted in the end of XIX century.
In spite of numerous achievements over the last more than hundred years on the classification of the equations of particular limited order (mainly up to fourth), the general problem for the higher order equations remains unsolved. This paper contributes towards solution of the arbitrary-order problem by building the complete Painlev\'e classification for the broad class of equations with the only constraint on the equation's degree but free of any order limitations.

The Painlev\'e classification of the non-linear polynomial ordinary differential equations
\begin{equation}
w^{(n)}=P(w^{(n-1)},w^{(n-2)},...,w,z),\label{pol}
\end{equation}
where $P$ is a polynomial in $w$ and its derivatives with coefficients locally
analytic in $z$, is known for the order $n\leq 4$.
For the order $n=1$ the well-known necessary and sufficient
condition of the Painlev\'e property for the equation (\ref{pol})
is $deg P=2$ (see for instance \cite{ince1956}, chapter XIII). For the
order $n=2$ classification has been built in the classical
works of Painlev\'e and Gambier (see for instance \cite{ince1956},
chapter XIV). For the order $n=3$ classification has been
started in the famous work of Chazy \cite{chazy1911equations} and recently
completed by C.Cosgrove \cite{cosgrove2000chazy}. Finally for the order $n=4$ the problem was solved by C.Cosgrove \cite{cosgrove2000P2, cosgrove2006P1}.
And although complete Painlev\'e classification has been successfully completed for certain algebraic classes of equations of the arbitrary order, such as binomial-type equations \cite{cosgrove1993binom, sobolevsky2005binomial3, sobolevsky2006binomialN, sobolevsky2006mono} and arbitrary algebraic equations that do not depend on the derivatives of order $n-1$ and $n-2$ \cite{sobolevskii2003algsing,sobolevsky2004stam,sobolevskii2005alg,sobolevsky2006mono}, the classification of analytically more simple polynomial-class ordinary differential equations of order $n\geq 5$ is not yet accomplished in the general case. 

Here we consider a class of polynomial ordinary differential
equations of arbitrary order $n\geq 2$, but restricted by the
degree of the right-hand side. Let $P$ be a quadratic form in $w$
and its derivatives, i.e. consider the equation
\begin{equation}
w^{(n)}=\sum\limits_{n-1\geq k\geq j\geq
0}a_{k,j}(z)w^{(k)}w^{(j)}+\sum\limits_{j=0}^{n-1}b_j(z)w^{(j)}+c(z),
\label{eqmain}
\end{equation}
where $a_{k,j},b_j,c$ are functions in $z$ analytic in the certain
complex domain $U$. Without loss of generality assume $a_{k,j}(z)\equiv a_{j,k}(z)$. We exclude the order $n=1$ since all Ricatti
equations, as mentioned afore, certainly possess the Painlev\'e
property.

First give the rigorous definition of the main concept -
the Painlev\'e property - for the considered equations (\ref{pol}) inline with  
it's classical understanding \cite{ince1956,conte1999painleve} but distinguishing commonly mixed concepts of freedom from movable branch points 
and non-polar singularities.

\textbf{Definition\,1.} \begin{it} Consider an arbitrary solution
$w=w(z)$ of the equation (\ref{eqmain}) analytic in the
neighborhood of some point $z^*\in U$ and a path $\Gamma$ with the
beginning in $z^*$ along which all the coefficients of (\ref{eqmain})
can be analytically continued while the analytical continuation of
$w(z)$ comes to a singularity. Such a singularity of $w(z)$ is
called movable singularity of the considered solution $w=w(z)$.
\end{it}

\textbf{Definition\,2.}
\begin{it}
The equation (\ref{eqmain}) is called to possess the Painlev\'e
property if the equation's solutions are single-valued near all of their movable
singularities. The equation (\ref{eqmain}) is called to possess
the strong Painlev\'e property if all (if any) movable singularities of it's
solutions are poles.
\end{it}

The equations possessing the strong Painlev\'e property according
to the definition 2 have been also introduced in \cite{gordoa2003new} as
the equations of Painlv\'e-type.

Our goal in the present paper is to find all of the equations of class (\ref{eqmain})
possessing the strong Painlev\'e property.

\section{First necessary
condition}

Introduce a constant $B=n-\max\{k+j:a_{k,j}(z)\not\equiv 0\}$. This characteristic
indicates the possible order of a movable pole for equation's (\ref{eqmain}) solutions. According to the theorem 4
\cite{Sobolevsky2004poly} if the initial equation (\ref{eqmain}) possesses
the Painlev\'e property and moreover the strong Painlev\'e
property, the number $B$ could be only $1$ or $2$. The corresponding two
possible forms of the equation (\ref{eqmain}) with the Painlev\'e
property are:
\begin{equation}
w^{(n)}=\sum\limits_{k=[n/2]}^{n-1}a_{k,n-1-k}(z)w^{(k)}w^{(n-1-k)}+\sum\limits_{k+j<n-1}a_{k,j}(z)w^{(k)}w^{(j)}+\sum\limits_{j=0}^{n-1}b_j(z)w^{(j)}+c(z),
 \label{eqB1}
\end{equation}
\begin{equation}
w^{(n)}=\sum\limits_{k=[(n-1)/2]}^{n-2}a_{k,n-2-k}(z)w^{(k)}w^{(n-2-k)}+\sum\limits_{k+j<n-2}a_{k,j}(z)w^{(k)}w^{(j)}+\sum\limits_{j=0}^{n-1}b_j(z)w^{(j)}+c(z),
\label{eqB2}
\end{equation}
where at least one of the coefficients $a_{k,j}$, for which
$n-(k+j)=B$, is not identically equal to zero. Moreover with respect to the theorem 3 \cite{Sobolevsky2004poly} at least one of such coefficients having $k\geq n-2$ or $l\geq n-2$ is not identically equal to zero as the equation's learning terms should include $w^{(n-1)}$ or $w^{(n-2)}$.

\section{Improved resonance condition}

According to \cite{Sobolevsky2004poly, sobolevskii2006modification} the equation
(\ref{pol}) always admits solutions with movable singularities.
And if the equation possesses the strong Painlev\'e property these
singularities are poles. Below we construct one general necessary
condition that the equation should satisfy in order for that to hold.

First of all if the equation (\ref{pol}) admits a solution with movable
singularity in a certain point $z=z_0$ being a pole, their should exist 
a Laurent expansion of the form
\begin{equation}
w=\sum\limits_{j=0}^\infty q_j (z-z_0)^{j-p}, \label{laurent}
\end{equation}
converging in a certain deleted neighborhood of $z=z_0$,
where $p$ is an integer order of the pole, while $a_j$ are complex coefficients, provided that $q_0\neq 0$. 
Further we first reproduce the well-known technique of the Painlev\'e test
\cite{ablowitz1980connection} in a rigorous way and also prove one additional necessary condition applied to the equation's so called resonance numbers.

Call some of the terms of the equation (\ref{pol}) leading
with respect to $(p,z_0)$ if their coefficients are nonzero in
$z_0$ and these terms produce the maximum order of the pole in
$z=z_0$ after substitution of the expression (\ref{laurent}) in
the equation (\ref{pol}) compared to all other terms of the equation. Denote the set of those leading
terms by $L(p,z_0)$.

Denote by $S$ the set of points in the complex domain being a zero
for at least one of equation's (\ref{pol}) (which could be written in the form (\ref{eqmain})) coefficients. Note that the set of leading terms $L(p,z_0)$
does not depend on $z_0$ for $z_0\not\in S$ and denote it by
$\hat{L}(p)$. For $z_0\in S$
the set $L(p,z_0)$ is a subset of $\hat{L}(z_0)$ if at least one
of the coefficients of the terms from $\hat{L}(p)$ is nonzero. In
the opposite case $L(p,z_0)$ is the set of new terms producing
lower pole orders in $z_0$.

Substituting series (\ref{laurent}) with undefined coefficients
$q_j$ into the equation (\ref{pol}) or (\ref{eqmain}) one should obtain an
algebraic equation
\begin{equation}
\label{determining}
T(q_0)=0
\end{equation}
for determination of $q_0$ (so called determining equation) and
the following recurrent equations for determining each further $q_j$ through
the earlier defined $q_0,q_1,...,q_{j-1}$:
\begin{equation}
R(j,q_0)q_j=Q_j(q_0,q_1,\ldots,q_{j-1}), \label{ajform}
\end{equation}
where $T$, $R$ and $Q_j$ are polynomials in all variables.

The polynomials $T$ and $R$ depend only on the leading terms
$L(p,z_0)$ (or $\hat{L}(p)$ in case $z_0\not\in S$), while the polynomials $Q_j$ together depend on all
the terms of the equation (\ref{pol}). For $z_0\not\in S$ the coefficients of the polynomials $T$ and $R$ are
polynomial expressions in all the values of the coefficients of the
leading terms from $\hat{L}(p)$ in the point $z_0$, while the
coefficients of $Q_j$ are polynomial expressions in $z_0$ and the
values of coefficients of (\ref{eqmain}) and their derivatives in
the point $z_0$. That is why for $z_0\not\in S$ the determining
equation (\ref{determining}) and the equations (\ref{ajform}) have the same 
general form for for all $z_0$ while coefficients of polynomials
$T,R,Q_j$ are analytic in $z_0$.

Since $q_0\neq 0$ in (\ref{laurent}) the determining equation
should admit nonzero roots and so there should exist at least two
leading terms $L(p,z_0)$ for each $z_0$ being the position of a
movable pole of order $p$. Generally the determining equation
(\ref{determining}) may have several solutions $q_0=q_0(z_0)$
being analytic functions in $z_0$ for all $z_0\not\in S$ probably
except of a countable set of algebraic singularities.

Starting from the choice of one of the roots $q_0$ of the
determining equation, one could use the equations (\ref{ajform})
to find all of the other $q_j$ except of those who's indexes $j=r$ are
the roots of the equation
\begin{equation}\label{resonans}
R(r,q_0)=0,
\end{equation}
called resonance equation for the given choice of $q_0$. Its roots
are called resonance numbers (or simply resonanses) of the equation (\ref{pol}) for the
selected $q_0$. One of those resonances is always $r=-1$
\cite{ablowitz1980connection} and is called trivial. The rigorous proof of
this fact could be found in \cite{conte1999painleve} (page 126-127).

The degree of the polynomial $R$ in $r$ does not exceed the maximum
order of the derivative of $w$ contained among the leading terms
$L(p,z_0)$. So the number of resonances other than trivial $r=-1$
is no more than $n-1$, while the maximum number of $n-1$ is
achieved only if the major derivative $w^{(n)}$ is contained among
the leading terms $L(p,z_0)$.

Denote the set of resonances corresponding to the given set of
$p,z_0,q_0$ by $r(p,z_0,q_0)$. As we have shown this set is
completely defined by the choice of $p,z_0,q_0$. Call the set of
numbers $p,z_0,q_0$ the initial characteristics of the
movable pole of the solution $w=w(z)$. So one can say that the
certain movable pole of the certain solution possesses the set of
resonances $r(p,z_0,q_0)$ defined by the initial characteristics
$p,z_0,q_0$ of this movable pole.

In that terms the well-known resonance condition for the Painlev\'e property \cite{conte1999painleve} can be formulated as following: 
for each $z_0$ and for every possible pair of $p,q_0$ (i.e. every pair $p,q_0$
such that $L(p,z_0)$ contains more than one term, $q_0$ satisfies
the corresponding determining equation and solutions of the
equation (\ref{pol}) with leading asymptotic behavior $w\sim
q_0(z-z_0)^{-p}$ really exist) all of the roots of
resonance equations (\ref{resonans}) with positive real part should be real integer.
Further in \cite{conte1999painleve} (page 132) for instance it is stated that in fact all the roots of (\ref{resonans}) should be distinct integers.

The positive integer resonances are values of indexes $j$ for
which the value of the coefficient $q_j$ can be arbitrary. The
further necessary conditions for the Painlev\'e property could be
obtained by analyzing the equations (\ref{ajform}) for indexes $j$
being a positive integer resonance. Such equations take the form
$Q_j(q_0,q_1,\ldots,q_{j-1})=0$ and give additional conditions
which should be satisfied with respect of the earlier obtained
expressions for the coefficients $q_0,q_1,\ldots,q_j$. These
conditions are called resonance conditions.

There exist a number of
works devoted to the Painlev\'e analysis of the certain equations
possessing negative resonances \cite{conte1999painleve} (page 139-140),
\cite{zdunek2000, conte1993perturbative, fordy1991analysing, musette1995non},
however it is still not quite clear what information could
negative resonances give for the Painlev\'e test in general.

It is also mentioned in \cite{ablowitz1980connection} that if for every pair of $p,q_0$ resonance roots do not contain $n-1$ nonnegative distinct integers then (\ref{laurent}) 
does not represent a general solution. It hints that perhaps general solution has a more complex shape and Painlev\'e property does not hold in such a case, however rigorous proof of this fact is not provided as for the author's knowledge. Below we'll fill this formal gap.

For a movable pole of order $p$ in point $z_0$ with leading coefficient $q_0$ 
call the set of resonances $r(p,z_0,q_0)$ complete if, besides the trivial resonance $-1$, it consists of 
$n-1$ distinct non-negative integers.

\textbf{Theorem\,1.} \begin{it} If the equation (\ref{pol})
possesses the strong Painlev\'e property and admits a solution
with movable pole in a certain point $z=z_1$ then for any deleted
neighborhood $V$ of $z_1$ there exists a solution of the equation
(\ref{pol}) with a movable pole $z_0\in V$ possessing complete set
of resonances.
\end{it}

In other words the theorem 1 states that negative resonances other
than trivial could exist but not for all movable singularities of
solutions of the equation (\ref{pol}) with strong Painlev\'e
property. Moreover positions of movable poles of solutions can not
be isolated and near each location of a movable pole one could always find other movable poles of other solutions 
possessing complete set of resonances.

\textbf{\emph{Proof of the theorem 1.}} Suppose the opposite case
i.e. let $w=w_0(z)$ be a solution of the equation (\ref{pol})
with a movable pole in $z=z_1$ and suppose that for any deleted
neighborhood $V$ of $z_1$ all movable poles of solutions of the
equation (\ref{pol}) located in $V$ possess no more than $n-2$
distinct non-negative integer resonances.

Consider an arbitrary small closed contour $\Gamma$ such as all the coefficients of the equation $\ref{pol}$ are analytic in $G=\overline{int\Gamma}$ (closed interior of $\Gamma$) while $z=z_0$ is the only singularity of the solution $w=w_0(z)$ in $G$. Select an arbitrary point $z=z^*$ on $\Gamma$ and consider a local representation for the general solution of equation (\ref{pol}) as $w=\phi(z,\lambda)$, where $\lambda=(w^0,w^1,\ldots,w^{n-1})$ is the set of parameters for the corresponding Cauchy initial value problem $w^{(j)}(z^*)=w^j$, while function $\phi$ is analytic for $z\in \Gamma$ and $\lambda$ from a certain $n$-dimensional neighborhood $V$ of $\lambda_0=(w_0(z^*),w'_0(z^*),\ldots,w^{(n-1)}_0(z^*))$. This way $\phi(z,\lambda_0)=w_0(z)$. For $\lambda$ close enough $\lambda_0$ solutions $w=\phi(z,\lambda)$ should also possess movable singularities inside $\Gamma$ (being poles since the equation possess strong Painlev\'e property) as otherwise function $w_0(z)$ would be analytic inside $\Gamma$ as a limit of the sequence of analytic functions). Without loss of generality assume that it holds for all $\lambda\in V$, otherwise simply narrow $V$ accordingly.

Then $\phi$ could be represented in the following form
\begin{equation}
\phi(z,\lambda)=h(z,\lambda)+\varphi(z,\lambda),
\label{genSol}
\end{equation}
where $\varphi(z,\lambda)=\frac{1}{2\pi i}\oint\limits_\Gamma \frac{\phi(\zeta,\lambda)}{\zeta-z}d\zeta$. Function $\varphi(z,\lambda)$ by definition is analytic in $G\times V$, while 
$h(z,\lambda)=\phi(z,\lambda)-\varphi(z,\lambda)$ for each $\lambda\in V$ is analytic in $z$ on $\Gamma$ and outside it including infinity, having the same poles inside $\Gamma$ as $\phi(z,\lambda)$ has. This way function $h$ is a rational function in $z$ with coefficients analytic for $\lambda\in V$. 

Let $\tau(\lambda)$ be the number of poles of function $h(z,\lambda)$ in $z$. Function $\tau$ is limited on $V$ since $h$ is rational. Consider $\tau^*=\max\limits_{\lambda\in V}\tau (\lambda)$ and a certain $\lambda^*\in V$ such as $\tau(\lambda^*)=\tau^*$. Then one can show that each of the poles $z_1\in G$ of the solution $\phi(z,\lambda^*)$ should possess a complete set of resonances. 

Indeed there exist a contour $\Gamma_1$ surrounding $z_1$ such as for any $\lambda$ from a certain neighborhood $V_1$ of $\lambda^*$, $\phi(z,\lambda)$ posses just one single pole inside $\Gamma_1$ as otherwise $\tau^*<\max\limits_{\lambda\in V}\tau (\lambda)$ (since for any $\lambda$ close enough to $\lambda^*$ function $h$ should also have at least one pole near each of all other locations of poles of $h(z,\lambda^*)$). Now similar to (\ref{genSol}) represent the general solution of (\ref{pol}) as
\begin{equation}
\phi(z,\lambda)=h_1(z,\lambda)+\varphi_1(z,\lambda),
\label{genSol1}
\end{equation}
where $\varphi_1$ is analytic for $z$ on and inside $\Gamma_1$ and $\lambda\in V_1$, while $h_1$ is the rational function in $z$ with coefficients analytic in $\lambda$ having a single pole in $z$. Obviously (\ref{genSol1}) is a local representation of the general solution for the equation (\ref{pol}) in the form of Laurent series around a movable pole. While if the resonance set of the solution $\phi(z,\lambda^*)$ corresponding to the pole in $z_1$ is not be complete, having only $n_1<n-1$ district positive integer resonance numbers, it would mean that all the coefficients of (\ref{genSol1}) could be uniquely defined through $m+1<n$ arbitrary coefficients, i.e. (\ref{genSol1}) would not represent the general solution of the equation (\ref{pol}).
Obtained contradiction completes the proof of the theorem 1.

\section{The resonance equation}

Now apply the resonance condition of theorem 1 to the two possible
cases (\ref{eqB1}) and (\ref{eqB2}) of the equation
(\ref{eqmain}). First from the equation (\ref{determining}) find the possible major coefficient $q_0$ for
the movable pole in the point $z=z_0$ in the form $q_0=1/f(z_0)$,
where
$$
f(z)=-\sum\limits_{k=[n/2]}^{n-1}a_{k,n-1-k}(z)k!(n-1-k)!/n!
$$
for the equation (\ref{eqB1}) and
$$
f(z)=\sum\limits_{k=[(n-1)/2]}^{n-2}a_{k,n-2-k}(z)(k+1)!(n-1-k)!/(n+1)!
$$
for the equation (\ref{eqB2}). Then obtain the resonance equations (\ref{resonans})
for the cases (\ref{eqB1}) and (\ref{eqB2}) in the corresponding
forms:
\begin{equation}
\begin{array}{c}
0=R(r)=\\
\prod\limits_{t=0}^{n-1}(r-1-t)-\sum\limits_{k=[n/2]}^{n-1}\frac{a_{k,n-1-k}(z_0)}{f(z_0)}\left(
(-1)^{n-1-k}(n-1-k)!\prod\limits_{t=0}^{k-1}(r-1-t)+(-1)^k
k!\prod\limits_{t=0}^{n-2-k}(r-1-t) \right) 
\end{array}
\label{resB1}
\end{equation}
and
\begin{equation}
\begin{array}{l}
0=R(r)=\prod\limits_{t=0}^{n-1}(r-2-t)-\\\sum\limits_{k=[(n-1)/2]}^{n-2}\frac{a_{k,n-2-k}(z_0)}{f(z_0)}\left(
(-1)^{n-2-k}(n-1-k)!\prod\limits_{t=0}^{k-1}(r-2-t)+(-1)^k
(k+1)!\prod\limits_{t=0}^{n-3-k}(r-2-t) \right)
\end{array}\label{resB2}
\end{equation}
(here and further the product of empty set of terms, as well as
$0!$ are considered to be equal to $1$).

Suppose that the initial equation (\ref{eqmain}) has the
Painlev\'e property. Then, according to the theorem 1, for
everywhere dense in the complex space set of $z_0$ the
corresponding resonance equation $R(r)=0$ has $n-1$ different
positive integer roots in addition to the trivial root $r=-1$. Since the
coefficients of the equation $R(r)=0$ as well as it's roots
depend continuously on $z_0$, the above mentioned condition is
possible only if these coefficients and roots are constant with
respect to $z_0$.

Denote the positive integer roots of the corresponding resonance
equation (\ref{resB1}) or (\ref{resB2}) by
$0<r_1<r_2<\ldots<r_{n-1}$ and let $r_0=-1$ be the trivial root.

\section{Solving the resonance equation in case of Bureau number $2$}

According to the Viet theorem for (\ref{resB2}) one can get
\begin{equation}
\sum\limits_{j=0}^{n-1}(r_j-2)=\sum\limits_{j=0}^{n-1}j,
\label{condB2-1}
\end{equation}
\begin{equation}
\sum\limits_{j=0}^{n-1}(r_j-2)^2=\sum\limits_{j=0}^{n-1}j^2+2h,
\label{condB2-2}
\end{equation}
where $h=a_{n-2,0}(z_0)/f(z_0)$ being constant with respect to
$z_0$. Of course $h$ should be integer. At the same time
$R(2)=(-1)^{n-1}(n-1)!h$. If $h>0$ then $R(2)$ and $R(-\infty)$
have different signs, so an interval $(-\infty;2)$ can not contain an even number of roots $r$ of the equation (\ref{resB2}). Since the only non-positive root is $r=-1$ this means that $r=1$ can not be another root. Also $R(2)\neq 0$, that is why
$2<r_1<r_2<\ldots<r_{n-1}$. Let $\delta_j=r_j-2-j$ for
$j=1,2,\ldots,n-1$. Then
$0\leq\delta_1\leq\delta_2\leq\ldots\leq\delta_{n-1}$ and, with
respect to (\ref{condB2-1}), one obtains
$0=\sum\limits_{j=0}^{n-1}(r_j-2-j)=-3+\sum\limits_{j=1}^{n-1}\delta_j$.
Then only the following three cases are possible:\\
1) $\delta_{n-3}=\delta_{n-2}=\delta_{n-1}=1$,
$\delta_1=\delta_2=\ldots=\delta_{n-4}=0$;\\
2) $\delta_{n-2}=1,\ \delta_{n-1}=2$,
$\delta_1=\delta_2=\ldots=\delta_{n-3}=0$;\\
3) $\delta_{n-1}=3$, $\delta_1=\delta_2=\ldots=\delta_{n-2}=0$.

However, with respect to (\ref{condB2-2}) one can obtain
\begin{equation}
h=\left(9+\sum\limits_{j=1}^{n-1}\left(\delta_j+j\right)^2-\sum\limits_{j=1}^{n-1}j^2\right)/2=\left(9+\sum\limits_{j=1}^{n-1}\delta_j^2+2\sum\limits_{j=1}^{n-1}j\delta_j\right)/2.
\label{contrB2-1}
\end{equation}
From the other hand
$(-1)^{n-1}(n-1)!h=R(2)=\left(2-(-1)\right)\prod\limits_{j=1}^{n-1}\left(2-(2+j+\delta_j)\right)=3(-1)^{n-1}\prod\limits_{j=1}^{n-1}\left(j+\delta_j\right)$,
consequently
\begin{equation}
h=3\prod\limits_{j=1}^{n-1}\left(j+\delta_j\right)/(n-1)!
\label{contrB2-2}
\end{equation}

In the case 1) with respect to (\ref{contrB2-1}) obtain
$h=\left(9+3+2(3n-6)\right)/2=3n$. From the other hand, according
to (\ref{contrB2-1}) obtain $h=3n/(n-3)$. It is possible only if
$n=4$. Then the corresponding initial equation (\ref{eqB2}) takes
form
\begin{equation}
w^{(IV)}=A(z)\left(ww''+(w')^2\right)+a_1(z)w'''+a_2(z)w''+a_3(z)w'w+a_4(z)w'+a_5(z)w^2+a_6(z)w+a_7(z).
\label{eqB2-1}
\end{equation}

In the case 2) with respect to (\ref{contrB2-1}) obtain
$h=\left(9+5+2\left(2(n-1)+(n-2)\right)\right)/2=3n+3$. From the
other hand, according to (\ref{contrB2-1}) obtain
$h=3(n+1)/(n-2)$. It is possible only if $n=3$. Then the
corresponding initial equation (\ref{eqB2}) takes form
\begin{equation}
w'''=A(z)ww'+a_1(z)w''+a_2(z)w'+a_3(z)w^2+a_4(z)w+a_5(z).
\label{eqB2-2}
\end{equation}

In the case 3) with respect to (\ref{contrB2-1}) obtain
$h=\left(9+9+2\left(3(n-1)\right)\right)/2=6+3n$. From the other
hand, according to (\ref{contrB2-1}) obtain $h=3(n+2)/(n-1)$. It
is possible only if $n=2$. Then the corresponding initial equation
(\ref{eqB2}) takes form
\begin{equation}
w''=A(z)w^2+a_1(z)w'+a_2(z)w+a_3(z). \label{eqB2-3}
\end{equation}

If $h<0$ then assume $\delta_j=r_j-2-j$ for $j=0,1,...,n-1$ and
according to (\ref{condB2-1}), (\ref{condB2-2}) obtain
$\sum\limits_{j=0}^{n-1}\delta_j=0$ and
$\sum\limits_{j=0}^{n-1}\left(\delta_j\right)^2+2\sum\limits_{j=0}^{n-1}\left(j\delta_j\right)=2h<0$.
However since $\delta_0\leq\delta_1\leq\ldots\leq\delta_{n-1}$ one
can obtain $\sum\limits_{j=0}^{n-1}\left(j\delta_j\right)>0$
because
$\sum\limits_{j=0}^{n-1}\left(j\delta_j\right)>\sum\limits_{j=0}^{n-1}\left(j\delta_{n-1-j}\right)$,
while
$\sum\limits_{j=0}^{n-1}\left(j\delta_j\right)+\sum\limits_{j=0}^{n-1}\left(j\delta_{n-1-j}\right)=(n-1)\sum\limits_{j=0}^{n-1}\delta_j=0$.
Consequently obtain the contradiction
$\sum\limits_{j=0}^{n-1}\left(\delta_j\right)^2=2h-2\sum\limits_{j=0}^{n-1}j\delta_j<0$.

Finally $h=0$ is not possible as $a_{n-2,0}(z)\not\equiv 0$. This way in case of Bureau number $2$ the only 3 possible forms for the equation (\ref{eqB2}) with the strong Painlev\'e property are (\ref{eqB2-1})-(\ref{eqB2-3}).

\section{Solving the resonance equation in case of Bureau number $1$}

First note that in case $a_{n-1,0}=0$ the equation (\ref{resB1})
after normalization by $r-1$ takes the form of (\ref{resB2}) with order
smaller by one and can be considered in the same way as above, which leads to the similar cases
1)-3) as for the equation (\ref{resB2}). For these cases the initial
equation (\ref{eqB1}) takes one of the  three corresponding forms:
\begin{equation}
w^{(V)}=A(z)\left(w'w'''+(w'')^2\right)+a_1(z)w^{(IV)}+S(w''',w'',w',w,z),
\label{eqB21-1}
\end{equation}
\begin{equation}
\begin{array}{c}
w^{(IV)}=A(z)w'w''+a_1(z)w'''+\left(a_2(z)w+a_3(z)\right)w''+a_4(z)\left(w'\right)^2+\\+a_5(z)ww'+a_6(z)w'+a_7(z)w^2+a_8(z)w+a_9(z),
\end{array}
\label{eqB21-2}
\end{equation}
\begin{equation}
w'''=A(z)\left(w'\right)^2+a_1(z)w''+a_2(z)ww'+a_3(z)w'+a_4(z)w^2+a_5(z)w+a_6(z),
\label{eqB21-3}
\end{equation}
where $S$ is a quadratic form in $w$ and its derivatives with
coefficients locally analytic in $z$ and cumulative number of
derivatives in each term not exceeding $3$.

Consider the case $a_{n-1,0}\neq 0$ in assumption that the
equation (\ref{resB1}) possesses $n-1$ distinct positive integer
roots being a necessary condition for the strong Painlev\'e property. Let those roots be $0<r_1<r_2<\ldots<r_{n-1}$ and also denote
$r_0=-1$. From the Viet theorem obtain
\begin{equation}
\sum\limits_{j=0}^{n-1}r_j=\sum\limits_{j=0}^{n-1}j+h,
\label{condB1-1}
\end{equation}
where $h=a_{n-1,0}(z_0)/f(z_0)$ being an integer non-zero constant. One can find
$R(1)=(-1)^n(n-1)!h\neq 0$. Since the interval $(-\infty;1)$ contains a
single non-multiple root $r=-1$ of the equation (\ref{resB1}), then $R(1)$ and $R(-\infty)$ should have different
signs, consequently $h<0$. Let $\delta_j=r_j-1-j$ for
$j=1,2,...,n-1$. Then $0\leq\delta_1\leq\ldots\leq\delta_{n-1}$
and according to (\ref{condB1-1}) obtain $-1\geq
h=r_0-1+\sum\limits_{j=1}^{n-1}(r_j-1-j)=-2+\sum\limits_{j=1}^{n-1}\delta_j$.
So only following two cases are possible:\\
1) $\delta_{n-1}=1$, $\delta_1=\delta_2=\ldots=\delta_{n-2}=0$,\\
2) $\delta_1=\delta_2=\ldots=\delta_{n-1}=0$.

However from one side $h=-2+\sum\limits_{j=1}^{n-1}\delta_j$ and
from the other
$(-1)^n(n-1)!h=R(1)=2(-1)^{n-1}\prod\limits_{j=1}^{n-1}(j+\delta_j)$.
In the case 1) one can find $h=-1$,
$(-1)^{n-1}(n-1)!=2(-1)^{n-1}\prod\limits_{j=1}^{n-1}(j+\delta_j)=2(-1)^{n-1}(n-2)!n$
and consequently obtain a contradiction $n-1=2n$.

In the case 2) the conditions
$h=-2+\sum\limits_{j=1}^{n-1}\delta_j$ and $(-1)^n(n-1)!h=R(1)$
correspond each other. The initial equation (\ref{eqB2}) takes the
following form
\begin{equation}
w^{(n)}=A(z)\left(w^2\right)^{(n-1)}+\sum\limits_{j+k<n-1}a_{k,j}(z)w^{(k)}w^{(j)}+\sum\limits_{j=0}^{n-1}b_j(z)w^{(j)}+c(z).
\label{eqB2-0}
\end{equation}

\section{The classification}

So we have proven the following

\textbf{Theorem\,2.} \begin{it} If the equation (\ref{eqmain})
possesses the strong Painlev\'e property it should be of one of
the following seven forms: (\ref{eqB2-1})-(\ref{eqB21-3}),
(\ref{eqB2-0}).
\end{it}

Note that only one of these seven aforementioned equations ---
equation (\ref{eqB2-0}) --- admits higher orders $n\geq 6$.

\textbf{Corollary\,1.} \begin{it} If the equation (\ref{eqmain})
of order $n\geq 6$ possesses the strong Painlev\'e property it
should be of the form (\ref{eqB2-0}).
\end{it}

\section{Necessary and sufficient conditions}

The theorem 2 gives only the necessary conditions for
the equation (\ref{eqmain}) to possess strong Painlev\'e
property. To complete the Painlev\'e classification one should
find the necessary and sufficient conditions of the strong
Painlev\'e property for each of the equations
(\ref{eqB2-1})-(\ref{eqB21-3}), (\ref{eqB2-0}).

Note that the equations (\ref{eqB2-1})-(\ref{eqB2-3}),
(\ref{eqB21-2}), (\ref{eqB21-3}) are already studied and the
conditions of Painlev\'e property for them are known. The equation
(\ref{eqB2-1}) corresponds to the class F-I \cite{cosgrove2000P2}, the
equation (\ref{eqB2-2}) --- to the Chazy Class XIII \cite{chazy1911equations},
\cite{cosgrove2000chazy}, the equation (\ref{eqB2-3}) --- to the
well-known case I(a) of the second order (\cite{ince1956}, chapter
XIV), equation (\ref{eqB21-2}) --- to the class F-VII
\cite{cosgrove2006P1} and the equation (\ref{eqB21-3})
--- to the Chazy Class I \cite{chazy1911equations}, \cite{cosgrove2000chazy}.

Consider the equation (\ref{eqB21-1}). Note that by means of
a variable substitution $w=-12v/A(z)$,
$u=v'''+6(v')^2-K_1(z)v''-K_2(z)v'v-K_3(z)v^2-K_4(z)v'-K_5(z)v$
for the certain choice of $K_1,K_2,...,K_5$ the equation
(\ref{eqB21-1}) can always be transformed to a system of the form
\begin{equation}
\left\{\begin{array}{l}
u''=L_1(z)u'+L_2(z)u+L_3(z)+h_1(z)vv'+h_2(z)v'+h_3(z)v\\
v'''=6(v')^2+K_1(z)v''+K_2(z)v'v+K_3(z)v^2+K_4(z)v'+K_5(z)v+u,\end{array}\right.
\label{eqB21-1s}
\end{equation}
where $L_1,L_2,L_3,h_1,h_2,h_3$ are locally analytic functions in
$z$.

\textbf{Lemma 1.} \begin{it} If the equation (\ref{eqB21-1})
possesses the Painlev\'e property, then $h_1(z)\equiv h_2(z)\equiv
h_3(z)\equiv 0$.
\end{it}

Indeed, introduce the small parameter transform $u=\alpha^{-4}U$,
$v=\alpha^{-1}V$, $z=z_0+\alpha x$, where $z_0$ is an arbitrary
constant, and obtain the transformed system in the form
\begin{equation}
\left\{\begin{array}{c} U''=\alpha L_1(z_0+\alpha x)U'+\alpha^2
L_2(z_0+\alpha x)U+\alpha^3 h_1(z_0+\alpha x)VV'+\alpha^4
h_2(z_0+\alpha x)V'+\\+\alpha^5 h_3(z_0+\alpha x)V +O(\alpha^6),\\
V'''=-6(V')^2+U+O(\alpha).\end{array}\right. \label{eqB21-1st}
\end{equation}
Consider the small parameter expansion for the solution of the
system (\ref{eqB21-1st}):
$$
\left\{\begin{array}{l} U=U_0(x)+\alpha U_1(x)+\alpha^2
U_2(x)+\alpha^3 U_3(x)+\alpha^4 U_4(x)+O(\alpha^5),
\\
V=V_0(x)+O(\alpha),
\end{array}
\right.
$$
where $U_0(x)$ is an arbitrary linear function, while $V=V_0(x)$
is an arbitrary solution of the equation
\begin{equation}
V'''=-6(V')^2+U_0(x), \label{eqB21-eqV0}
\end{equation}
Then $U_1$ and $U_2$ can be found as polynomials in $x$, while
$U_3''(x)=H_3(x)+h_1(z_0)V_0(x)V_0'(x)$, where $H_3$ is a
polynomial in $x$.

So if the equation (\ref{eqB21-1}) possesses the Painlev\'e
property then either $h_1(z_0)=0$, or for any solution $V=V_0(x)$
of the equation (\ref{eqB21-eqV0}) the expression $\int\int
V_0(x)V_0'(x)dxdx$ should be single-valued so the function
$V_0(x)^2$ should always possess zero residue in any of its
singularities. However, analyzing the possible Laurent series
representation of general solution of (\ref{eqB21-eqV0}) near
movable pole one can find that this suggestion is invalid. So
$h_1(z_0)=0$ for arbitrary $z_0$, consequently $h_1(z)\equiv 0$.

Further, in case $h_1(z)\equiv 0$, one can find
$U_4''(x)=H_4(x)+h_2(z_0)V_0'(x)$ where $H_4$ is a polynomial in
$x$. So if $h_2(z_0)\neq 0$, the expression $\int V_0(x)dx$ should
be single-valued for any solution $V=V_0(x)$ of the equation
(\ref{eqB21-eqV0}). However this expression is multi-valued near
the first order movable poles of $V_0$. Consequently $h_2(z)\equiv
0$. Finally $h_3(z)\equiv 0$ in the same way, since
$U_5''(x)=H_5(x)+h_3(z_0)V_0(x)$ while the expression $\int\int
V_0(x)dxdx$ is also multi-valued near the first order movable
poles of $V_0$. The proof of lemma 1 is now complete.

Now from lemma 1 one can see that the equation (\ref{eqB21-1})
with the Painlev\'e property should necessary possesses the second
integral
\begin{equation}
v'''=-6(v')^2+K_1(z)v''+K_2(z)v'v+K_3(z)v^2+K_4(z)v'+K_5(z)v+u(z),
\label{eqB21-1_int}
\end{equation}
where $u(z)$ contains two constants of integration being an
arbitrary solution of the second order linear equation
\begin{equation}
u''=L_1(z)u'+L_2(z)u+L_3(z). \label{eqB21-1_lin}
\end{equation}
The integral (\ref{eqB21-1_int}) is the Chazy Class I equation
\cite{chazy1911equations}, \cite{cosgrove2000chazy} and the necessary and sufficient
conditions of the strong Painlev\'e for it are well-known:
$K_1(z)\equiv K_2(z)\equiv 0$, $K_3(z)\equiv K_4(z)$,
$K_4''(z)\equiv \left(K_4(z)\right)^2$, $K_5''(z)\equiv
K_4(z)K_5(z)$ and $u''(z)\equiv K_4(z)u(z)/3+
\left(K_5(z)/6\right)^2$, i.e. $L_1(z)\equiv 0$, $L_2(z)\equiv
K_4(z)/3$, $L_3(z)\equiv \left(K_5(z)/6\right)^2$.

Finally consider the equation (\ref{eqB2-0}).

\textbf{Theorem\,3.}
\begin{it}
The equation (\ref{eqB2-0}) possess the strong Painlev\'e property
if and only if it is linearizable by means of the variable change
\begin{equation}
u=w'-A(z)w^2-B(z)w, \label{eqB2-0_varchange}
\end{equation}
where $B(z)$ is a certain locally analytic function.
\end{it}

Of course if the equation (\ref{eqB2-0}) is linearizable by means
of (\ref{eqB2-0_varchange}) then the equation (\ref{eqB2-0})
possess the strong Painlev\'e property, since this way the general solution $u$ of the linear differential equation is free of any movable singularities, 
while the correspondent function $w$ can be found by resolving the Ricatti equation
(\ref{eqB2-0_varchange}) and so all of it's movable singularities are poles.

To prove the inverse statement consider the equation
(\ref{eqB2-0}) and assume that it possesses the strong Painlev\'e
property and consequently the Painlev\'e property. Introduce the
variable change (\ref{eqB2-0_varchange}) where $B(z)$ is a locally
analytic function, undefined yet. Then the equation (\ref{eqB2-0})
can be transformed to:
\begin{equation}
\label{eqB2-0_s1}
\begin{array}{l}
u^{(n-1)}=\sum\limits_{p(\chi)\leq
n}\tilde{a}_{\chi}(z)w^{\chi_0}\prod\limits_{j=1}^{n-1}\left(u^{(j-1)}\right)^{\chi_j},\\
w'=u+A(z)w^2+B(z)w,
\end{array}
\end{equation}
where $\chi=(\chi_0,\chi_1,...,\chi_{n-1})$ are multi-indexes with
integer non-negative components,
$p(\chi)=\sum\limits_{j=0}^{n-1}(j+1)\chi_j$, and $\tilde{a}_\chi$
are locally analytic in $z$ coefficients which can be polynomially
expressed in terms of the coefficients of the initial equation
(\ref{eqB2-0}), functions $A,B$ and their derivatives. Consider
the coefficient $\tilde{a}_{(0,0,...,0,n)}$ at $w^n$ in the
right-hand side of the first equation of the system
(\ref{eqB2-0_s1}). It depends linearly on $B(z)$ being of the form
$-n!B(z)A(z)^n+...$, where dots denote terms not containing
$B(z)$. Consequently by the corresponding choice of $B(z)$, one
can always make the degree of the system (\ref{eqB2-0_s1}) first equation's right-hand side to be not higher than
$n-1$ with respect to $w$.

Demonstrate that in this case the right-hand side of the first equation of
the system (\ref{eqB2-0_s1}) should not depend on $w$ at all.
Indeed, suppose the opposite case. Then by means of a small parameter
transform $z=z_0+\alpha x$, $u=\alpha^{-2}U$, $w=\alpha^{-1}W$
with arbitrary constant $z_0$, the system (\ref{eqB2-0_s1}) is
transformed to
\begin{equation}
\left\{\begin{array}{l} \frac{d^{n-1}U}{dx^{n-1}}=\alpha
T(U^{(n-2)},U^{(n-3)},...,U,x,\alpha)+\alpha^k
H(W,U^{(n-2)},U^{(n-3)},...,U)+o(\alpha^k),
\\[0.5ex]
\frac{dW}{dx}=U+A(z_0)W^2+O(\alpha),
\end{array}
\right. \label{eqB2-0_s2}
\end{equation}
where $k$ is a certain positive integer, while $T,H$ are
polynomials in all variables, while $1\leq {\rm deg}_W H\leq n-1$.
The solution of the system (\ref{eqB2-0_s2}) can be found in the
form
\begin{equation}
\label{eqB2-0_sol}
\begin{array}{c}
U(x,\alpha)=U_0(x)+\sum\limits_{j=1}^{k-1}U_j(x)\alpha^j+\alpha^k\int\int\ldots\int
H(W_0(x),U_0^{(n-2)}(x),U_0^{(n-3)}(x),...,U_0(x))dx^{n-1}+\\+o(\alpha^k),\\[0.5ex]
W(x,\alpha)=W_0(x)+O(\alpha),
\end{array}
\end{equation}
where $U_0(x)=C_{n-2}x^{n-2}+C_{n-3}x^{n-3}+...+C_1 x+C_0$ and
$W_0$ is an arbitrary solution of the Ricatti equation
$W'=A(z_0)W^2+U_0(x)$ with a movable simple pole in $x=-C$, while
$C,C_0,C_1,...,C_{n-2}$ --- are arbitrary complex constants. One
can see that the function $U(x,\alpha)$, being determined by
(\ref{eqB2-0_sol}) is multi-valued in some neighborhood of the
point $x=-C$ for sufficiently close to zero nonzero $\alpha$,
because the function
$H(W_0(x),U_0^{(n-2)}(x),U_0^{(n-3)}(x),...,U_0(x))$ in general
case admits a pole of order ${\rm deg}_W H\leq n-1$ in $x=-C$.

Consequently (\ref{eqB2-0_varchange}) transforms the equation
(\ref{eqB2-0}) to the polynomial differential equation of order
$n-1$ in $u$, not depending on $w$. In case this equation is nonlinear, then according to the theorem 4
\cite{Sobolevsky2004poly} its Bureau number should be $1$ or $2$ so for at
least one of the terms the inequality $p(\chi)\geq n+1$ should
hold. However for all terms we have $p(\chi)\leq n$. This
contradiction completes the proof of the theorem 3.

This way the necessary and sufficient conditions of the strong
Painlev\'e property for each of the possible seven cases
(\ref{eqB2-1})-(\ref{eqB21-3}), (\ref{eqB2-0}) are constructed completing the 
strong-Painlev\'e classification for the second
degree arbitrary order polynomial equations (\ref{eqmain}). Note that the only possible equation (\ref{eqmain}) of
order $n\geq 6$ with the strong Painlev\'e property, i.e. the
equation (\ref{eqB2-0}), is linearizable, while others could be transformed to the equations previously known. In particular this means that solutions of the
second degree polynomial differential equations having the strong Painlev\'e property do not provide any new transcendental functions.

\section*{Conclusions}
We've built a complete classification of the second degree arbitrary order polynomial ordinary differential equations (\ref{eqmain}) having strong Painlev\'e property. We proved that all such equations are contained in 7 classes (\ref{eqB2-1})-(\ref{eqB21-3}),
(\ref{eqB2-0}), and for each of those classes necessary and sufficient conditions for the strong Painlev\'e property are obtained. Six classes (\ref{eqB2-1})-(\ref{eqB21-3}) happen to have a limited order $n\leq 6$ and if having a strong Painlev\'e property all appear to be reduced to the previously known equations, while the only class (\ref{eqB2-0}) of unlimited order appears to be linearizable.
This way it is proven that second degree polynomial equations (\ref{eqmain}) of the arbitrary order having the strong Painlev\'e property are all integrable by means of known functions and do not provide any new transcendental solutions.

\bibliography{Painleve}

\begin{thebibliography}{10}
\providecommand{\url}[1]{\texttt{#1}}
\providecommand{\urlprefix}{URL }
\expandafter\ifx\csname urlstyle\endcsname\relax
  \providecommand{\doi}[1]{doi:\discretionary{}{}{}#1}\else
  \providecommand{\doi}{doi:\discretionary{}{}{}\begingroup
  \urlstyle{rm}\Url}\fi
\providecommand{\bibAnnoteFile}[1]{%
  \IfFileExists{#1}{\begin{quotation}\noindent\textsc{Key:} #1\\
  \textsc{Annotation:}\ \input{#1}\end{quotation}}{}}
\providecommand{\bibAnnote}[2]{%
  \begin{quotation}\noindent\textsc{Key:} #1\\
  \textsc{Annotation:}\ #2\end{quotation}}
\providecommand{\eprint}[2][]{\url{#2}}

\bibitem{ince1956}
Ince E (1956) \selectlanguage{russian} Ordinary differential equations.
\newblock Dover, New York.
\bibAnnoteFile{ince1956}

\bibitem{chazy1911equations}
Chazy J (1911) Sur les {\'e}quations diff{\'e}rentielles du troisi{\`e}me ordre
  et d'ordre sup{\'e}rieur dont l'int{\'e}grale g{\'e}n{\'e}rale a ses points
  critiques fixes.
\newblock Acta Mathematica 34: 317--385.
\bibAnnoteFile{chazy1911equations}

\bibitem{cosgrove2000chazy}
Cosgrove CM (2000) Chazy classes ix--xi of third-order differential equations.
\newblock Studies in applied mathematics 104: 171--228.
\bibAnnoteFile{cosgrove2000chazy}

\bibitem{cosgrove2000P2}
Cosgrove CM (2000) Higher-order painlev{\'e} equations in the polynomial class
  i. bureau symbol p2.
\newblock Studies in applied mathematics 104: 1--65.
\bibAnnoteFile{cosgrove2000P2}

\bibitem{cosgrove2006P1}
Cosgrove CM (2006) Higher-order painlev{\'e} equations in the polynomial class
  ii: Bureau symbol p1.
\newblock Studies in Applied Mathematics 116: 321--413.
\bibAnnoteFile{cosgrove2006P1}

\bibitem{cosgrove1993binom}
Cosgrove CM (1993) All-binomial-type painlev{\'e} equations of the second order
  and degree three or higher.
\newblock Studies in applied mathematics 90: 119--187.
\bibAnnoteFile{cosgrove1993binom}

\bibitem{sobolevsky2005binomial3}
Sobolevsky S (2005) Binomial-type ordinary differential equations of the third
  order.
\newblock Studies in Applied Mathematics 114: 1--15.
\bibAnnoteFile{sobolevsky2005binomial3}

\bibitem{sobolevsky2006binomialN}
Sobolevsky S (2006) Painlev{\'e} classification of binomial type ordinary
  differential equations of the arbitrary order.
\newblock Studies in Applied Mathematics 117: 215--237.
\bibAnnoteFile{sobolevsky2006binomialN}

\bibitem{sobolevsky2006mono}
Соболевский С (2006) Подвижные особые точки решений обыкновенных
  дифференциальных уравнений.
\newblock Минск: БГУ (in Russian), 118 pp.
\bibAnnoteFile{sobolevsky2006mono}

\bibitem{sobolevskii2003algsing}
Sobolevskii S (2003) Movable singular points of ordinary differential equations
  with algebraic singularities of the right-hand side.
\newblock Differential Equations 39: 381--386.
\bibAnnoteFile{sobolevskii2003algsing}

\bibitem{sobolevsky2004stam}
Sobolevsky S (2004) Movable singularities of a class of nonlinear ordinary
  differential equations of arbitrary order.
\newblock Studies in Applied Mathematics 112: 227--234.
\bibAnnoteFile{sobolevsky2004stam}

\bibitem{sobolevskii2005alg}
Sobolevskii S (2005) Movable singular points of algebraic ordinary differential
  equations.
\newblock Differential Equations 41: 1146--1154.
\bibAnnoteFile{sobolevskii2005alg}

\bibitem{conte1999painleve}
Conte R (1999) The Painlev{\'e} property: one century later, volume~1.
\newblock Springer New York.
\bibAnnoteFile{conte1999painleve}

\bibitem{gordoa2003new}
Gordoa PR, Joshi N, Pickering A (2003) A new technique in nonlinear singularity
  analysis.
\newblock Publications of the Research Institute for Mathematical Sciences 39:
  435--449.
\bibAnnoteFile{gordoa2003new}

\bibitem{Sobolevsky2004poly}
Sobolevskii S (2004) Movable singular points of polynomial ordinary
  differential equations.
\newblock Differential Equations 40: 807--814.
\bibAnnoteFile{Sobolevsky2004poly}

\bibitem{sobolevskii2006modification}
Sobolevskii S (2006) On a modification of the small parameter method.
\newblock Differential Equations 42: 218--228.
\bibAnnoteFile{sobolevskii2006modification}

\bibitem{ablowitz1980connection}
Ablowitz MJ, Ramani A, Segur H (1980) A connection between nonlinear evolution
  equations and ordinary differential equations of p-type. i.
\newblock Journal of Mathematical Physics 21: 715--721.
\bibAnnoteFile{ablowitz1980connection}

\bibitem{zdunek2000}
Здунек А, Мартынов И, Пронько В (2000) О рациональных решениях дифференциальных
  уравнений.
\newblock Вестник Гродненского Государственного университета им Янки Купалы Сер
  2 Математика (in Russian) 1: 33-39.
\bibAnnoteFile{zdunek2000}

\bibitem{conte1993perturbative}
Conte R, Fordy AP, Pickering A (1993) A perturbative painlev{\'e} approach to
  nonlinear differential equations.
\newblock Physica D: Nonlinear Phenomena 69: 33--58.
\bibAnnoteFile{conte1993perturbative}

\bibitem{fordy1991analysing}
Fordy A, Pickering A (1991) Analysing negative resonances in the painlev{\'e}
  test.
\newblock Physics Letters A 160: 347--354.
\bibAnnoteFile{fordy1991analysing}

\bibitem{musette1995non}
Musette M, Conte R (1995) Non-fuchsian extension to the painlev{\'e} test.
\newblock Physics Letters A 206: 340--346.
\bibAnnoteFile{musette1995non}

\end{thebibliography}

\end{document}